\input amstex
\documentstyle{amsppt}
\magnification=\magstep1
\pagewidth{15.9truecm}
%%\hcorrection{-0.25truein}
\pageheight{23.8truecm}
%%\vcorrection{0.75truein}

\input epsf.tex

%%\NoRunningHeads
%%\NoPageNumbers

\NoBlackBoxes
\nologo

\long\def\ignore#1\endignore{\par XYPIC DIAGRAM\par}
\long\def\ignore#1\endignore{#1}

\ignore
\input xy \xyoption{matrix} \xyoption{arrow}
%%          \xyoption{curve}
\def\edge{\ar@{-}}

\def\place{*{}+0}
\def\plb{*{\bullet}+0}
\def\plc{*{\circ}+0}

\def\cbb{\plc&\plb&\plb}
\def\bbb{\plb&\plb&\plb}
\def\phcirc{\phantom{\circ}}

\def\hhrz{\place \edge[rr] &\place &\place}

\def\hhrzvrt{\place \edge[rr] \edge[d] &\place \edge[d]
     &\place \edge[d]}

\endignore

\def\bfq{\bold q}
\def\Xbar{\overline{X}}
\def\CC{{\Bbb C}}
\def\QQ{{\Bbb Q}}
\def\ZZ{{\Bbb Z}}
\def\gfrak{{\frak g}}

\def\M{{\Cal M}}
\def\O{{\Cal O}}
\def\X{{\Cal X}}

\def\eqdef{\mathrel{\overset\roman{def}\to=}}
\def\Oq{{\Cal O}_q}
\def\Oqktwo{\Oq(k^2)}
\def\Oktwo{{\Cal O}(k^2)}
\def\Obfqkn{{\Cal O}_\bfq(k^n)}
\def\Obfqkxn{{\Cal O}_\bfq((\kx)^n)}
\def\Oqkn{\Oq(k^n)}
\def\Okn{{\Cal O}(k^n)}
\def\OqMtwo{\Oq(M_2(k))}
\def\OqSLtwo{\Oq(SL_2(k))}
\def\OqGLtwo{\Oq(GL_2(k))}
\def\OqMn{\Oq(M_n(k))}
\def\OqSLn{\Oq(SL_n(k))}
\def\OqGLn{\Oq(GL_n(k))}
\def\kbar{\overline{k}}
\def\kx{k^{\times}}
\def\kxtwo{(\kx)^2}
\def\kxfour{(\kx)^4}
\def\kxn{(\kx)^n}

\def\spec{\operatorname{spec}}
\def\Hspec{{H}\operatorname{-spec}}
\def\specj{\operatorname{spec}_J}
\def\prim{\operatorname{prim}}
\def\stabh{\operatorname{Stab}_H}
\def\E{{\Cal E}}
\def\Aut{\operatorname{Aut}}

\def\Fract{\operatorname{Fract}}
\def\End{\operatorname{End}}
\def\chr{\operatorname{char}}
\def\rank{\operatorname{rank}}

\def\Ami{{\bf 1}}
\def\BroGoo{{\bf 2}}
\def\BrGo{{\bf 3}}
\def\BrVe{{\bf 4}}
\def\Cau{{\bf 5}}
\def\ChPr{{\bf 6}}
\def\DeCEP{{\bf 7}}
\def\DKP{{\bf 8}}
\def\Dix{{\bf 9}}
\def\GMurcia{{\bf 10}}
\def\GDurham{{\bf 11}}
\def\GLenMurcia{{\bf 12}}
\def\GLenDuke{{\bf 13}}
\def\GLenIJM{{\bf 14}}
\def\GLenJalg{{\bf 15}}
\def\qaffine{{\bf 16}}
\def\specstrat{{\bf 17}}
\def\qaffquo{{\bf 18}}
\def\GoSt{{\bf 19}}
\def\Hod{{\bf 20}}
\def\HoLethree{{\bf 21}}
\def\HoLeN{{\bf 22}}
\def\Hortdiss{{\bf 23}}
\def\Hortpaper{{\bf 24}}
\def\Ing{{\bf 25}}
\def\IrSm{{\bf 26}}
\def\Jsur{{\bf 27}}
\def\Kas{{\bf 28}}
\def\KlSc{{\bf 29}}
\def\Lau{{\bf 30}}
\def\McRo{{\bf 31}}
\def\Moe{{\bf 32}}
\def\Mus{{\bf 33}}
\def\Nor{{\bf 34}}
\def\RTF{{\bf 35}}
\def\Van{{\bf 36}}

\topmatter

\title Quantized Coordinate Rings and Related Noetherian Algebras
\endtitle

\rightheadtext{Quantized Coordinate Rings}

\author K. R. Goodearl \endauthor

\address Department of Mathematics, University of California, Santa
Barbara, CA 93106, USA \endaddress

\email goodearl\@math.ucsb.edu \endemail

\abstract  This paper contains a survey of some ring-theoretic
aspects of quantized coordinate rings, with primary focus on the prime
and primitive spectra. For these algebras, the overall structure of
the prime spectrum is governed by a partition into strata determined
by the action of a suitable group of automorphisms of the algebra. We
discuss this stratification in detail, as well as its use in
determining the primitive spectrum -- under suitable conditions, the
primitive ideals are precisely those prime ideals which are maximal
within their strata. The discussion then turns to the global
structure of the primitive spectra of quantized coordinate rings, and
to the conjecture that these spectra are topological quotients of the
corresponding classical affine varieties. We describe the solution to
the conjecture for quantized coordinate rings of full affine spaces
and (somewhat more generally) affine toric varieties. The final part
of the paper is devoted to the quantized coordinate ring of $n\times
n$ matrices. We mention parallels between this algebra and the
classical coordinate ring, such as the primeness of quantum analogs
of determinantal ideals. Finally, we describe recent work which
determined, for the $3\times3$ case, all prime ideals invariant under
the group of winding automorphisms governing the stratification
mentioned above.
\endabstract

\thanks This is an expository paper, based on work published elsewhere.
\endthanks

\endtopmatter

\document

\head Introduction \endhead

First, a {\it caveat\/} concerning the title: This survey is not
designed to be either an introduction to or a discussion of quantum
groups. Rather, we present some of the ring theory that has arisen in
studying the structure of certain algebras found among quantum
groups. Here we only give a few words of background, and later we
present some representative examples. An introduction to the general
theory of quantum groups can be found in many books; as a
small sample, we mention \cite{\BrGo, \ChPr, \Kas, \KlSc}.

The term `quantized coordinate rings' refers to certain algebras
that, loosely speaking, are deformations of the classical coordinate
rings of affine algebraic varieties or algebraic groups. These algebras
are typically not commutative, but they turn out to have many other
properties analogous to the classical case -- for example, they are
noetherian, and most of the ones that have been introduced to date are
integral domains, with finite global dimension. To take the most basic
case, recall that the classical coordinate ring of affine $n$-space
over a field $k$ is just a polynomial ring in $n$ indeterminates over
$k$. Thus, a `quantized' coordinate ring of affine $n$-space should
be some type of noncommutative polynomial ring in $n$ indeterminates,
such as an $n$-fold iterated skew polynomial extension of $k$. For
the canonical examples, see Section 1.1.
\medskip

These notes are arranged in three parts, which focus on prime
ideals, primitive ideals, and the quantized coordinate rings of
matrices, respectively. Most of the material in Parts I and II is
excerpted from \cite{\BrGo}, where the reader can find a much more
detailed development. The aim of Part III is to illustrate how the
general picture developed in the first two parts applies to a
particularly interesting quantized coordinate ring; the discussion is
taken partly from
\cite{\BrGo} and partly from the recent paper \cite{\GLenJalg}. 

Throughout, we work over a base field $k$, and our parameters will be
elements of $\kx$, that is, nonzero scalars from $k$. The
characteristic of $k$ may be arbitrary, and for many results, it does
not matter whether or not $k$ is algebraically closed. We will
often concentrate on the so-called {\it generic\/} case, meaning that
our parameters are not roots of unity, but when not specified, the
parameters may be arbitrary. The key difference is that when
sufficiently many parameters are roots of unity, quantized coordinate
rings are finitely generated modules over their centers, and their
study proceeds via the theory of rings with polynomial identity. Our
aim here is to concentrate on the non-PI case, which requires very
different tools (some not yet invented).

\head I. Prime Ideals\endhead

In Part I, we concentrate on prime ideals in quantized
coordinate rings and related algebras, more precisely, on ways to
organize the {\it prime spectrum\/} -- the set $\spec A$ of prime
ideals in an algebra $A$. We view $\spec A$ not just as a set, but as
a topological space, equipped with the standard Zariski topology.

In order to have available a few examples with which to illustrate the
results and techniques, we begin by presenting some of the standard
quantized coordinate rings. For a survey of most of the known types,
see \cite{\GMurcia}.

\definition{1.1\. Some quantized coordinate rings} Let $q\in\kx$. The
{\it quantized coordinate ring of the $xy$-plane\/} with parameter $q$
is the
$k$-algebra
$$\Oqktwo \eqdef k\langle x,y \mid xy=qyx\rangle.$$
This algebra is often called a {\it quantum plane\/} for short.

Quantized coordinate rings for higher-dimensional spaces are defined
similarly, except that more choices of parameters are allowed. Let
$\bfq= (q_{ij})$ be a {\it multiplicatively
antisymmetric\/} $n\times n$ matrix over $k$, meaning that
$q_{ii}=1$ and
$q_{ji}= q_{ij}^{-1}$ for all $i,j$. The {\it quantized coordinate
ring of affine $n$-space\/} with parameter matrix $\bfq$ is the
$k$-algebra
$$\Obfqkn \eqdef k\langle x_1,\dots,x_n \mid x_ix_j= q_{ij}x_jx_i
\text{\ for all\ } i,j\rangle.$$
There is a single-parameter version of this algebra, defined for
$q\in\kx$ as follows:
$$\Oqkn \eqdef k\langle x_1,\dots,x_n \mid x_ix_j= qx_jx_i
\text{\ for all\ } i<j\rangle.$$
This is the special case of $\Obfqkn$ for which the matrix $\bfq$ has
the form 
$$\left(\smallmatrix 1&q&q&\cdots&q&q\\ q^{-1}&1&q&\cdots&q&q\\
q^{-1}&q^{-1}&1&\cdots&q&q\\ &\vdots&&&\vdots\\
q^{-1}&q^{-1}&q^{-1}&\cdots&1&q\\ 
q^{-1}&q^{-1}&q^{-1}&\cdots&q^{-1}&1 \endsmallmatrix\right).$$

The {\it quantized coordinate ring of $2\times2$ matrices\/} with
parameter $q$ is the
$k$-algebra $\OqMtwo$ given by four generators $X_{11}$, $X_{12}$,
$X_{21}$,
$X_{22}$ and six relations
$$\xalignat2 X_{11}X_{12} &= qX_{12}X_{11} &X_{12}X_{22} &=
qX_{22}X_{12} \\
X_{11}X_{21} &= qX_{21}X_{11} &X_{21}X_{22} &= qX_{22}X_{21} \\
X_{12}X_{21} &= X_{21}X_{12} 
&X_{11}X_{22}-X_{22}X_{11} &= (q-q^{-1})X_{12}X_{21}. \endxalignat$$
The first five relations, which are all of the form $xy=ryx$ for
generators $x$ and $y$ and scalars
$r$, can be summarized in the following mnemonic diagram:
\ignore
$$\xymatrix{
X_{11} \ar[r]^q \ar[d]_q &X_{12} \ar[d]^q \ar@{<->}[dl]_1 \\
X_{21}\ar[r]^q &X_{22}
}$$
\endignore

The element $D_q \eqdef X_{11}X_{22}-qX_{12}X_{21}$ in $\OqMtwo$ is
called the ($2\times2$) {\it quantum determinant\/}; as is easily
checked, $D_q$ lies in the center of $\OqMtwo$. The {\it quantized
coordinate rings of $GL_2(k)$ and $SL_2(k)$\/} are the algebras
$$\OqGLtwo \eqdef \OqMtwo[ D_q^{-1}] \qquad \text{and} \qquad
\OqSLtwo \eqdef \OqMtwo/
\langle D_q-1 \rangle.$$
Analogous algebras $\OqMn$, $\OqSLn$, and $\OqGLn$ have been defined
for arbitrary $n$, but we shall not give their definitions until
later (Section 3.1). 
\enddefinition

A general principle from the study of quantum phenomena in physics,
which holds equally well in mathematical studies of quantum algebras,
is that {\it quantization destroys symmetry\/}, meaning that a
quantized version of a classical system (physical or mathematical)
tends to be more rigid, with less symmetry. We illustrate this
principle with coordinate rings of the plane $k^2$. For instance, the
classical coordinate ring $k[x,y]$ has a huge supply of prime ideals,
but the quantized coordinate ring has far fewer, as the following
example shows. The same can also be said for automorphisms, as we
shall see shortly.

\definition{1.2\. Example}
When $k$ is algebraically closed and $q$ is not a root of unity, the
prime spectrum of $\Oqktwo$ can be displayed as follows:

\ignore
$$\xymatrixrowsep{4pc}\xymatrixcolsep{2pc}
\xymatrix{
 \ar@{--}[r] &\langle x,\, y-\beta \rangle \save+<0ex,-4ex>
\drop{(\beta \in k^\times)} \restore \ar@{--}[r] &&\langle x,y \rangle
&\ar@{--}[r] &\langle x-\alpha,\, y \rangle
\save+<0ex,-4ex> \drop{(\alpha \in k^\times)} \restore
\ar@{--}[r] & \\
 &\langle x\rangle \edge[ul] \edge[ur] \edge[urr] &&&&\langle y\rangle
\edge[ull] \edge[ul] \edge[ur]\\
 &&&\langle 0\rangle \edge[ull] \edge[urr]
}$$
\endignore
\enddefinition

Another difference in symmetry between the classical and quantized
coordinate rings of the plane is found in the automorphisms of these
algebras. As an algebraic variety, the plane is completely homogeneous,
in that any point can be moved to any other point by a translation.
These translations induce automorphisms of $k[x,y]$ of the following
form: For any scalars $a,b\in k$, there is a $k$-algebra
automorphism of $k[x,y]$ such that $x\mapsto x+a$ and $y\mapsto y+b$.
In the quantum case, however, $\Oqktwo$ has no such automorphisms
except the identity (corresponding to $a=b=0$). Fortunately, $\Oqktwo$
is not bereft of automorphisms -- there are multiplicative analogs of
the translation automorphisms, mapping $x$ and $y$ to scalar multiples
of themselves. In fact, all of our standard examples have a supply of
automorphisms of this type, as follows.

\definition{1.3\. Some automorphisms} We define some families of
$k$-algebra automorphisms on the quantized coordinate rings discussed
above. Each family of automorphisms is parametrized by tuples of
nonzero scalars, i.e., by elements from one of the multiplicative
groups $(\kx)^r$.

For $(\alpha,\beta)\in \kxtwo$, 
there is an automorphism of $\Oqktwo$ such that
$x\mapsto \alpha x$ and $y\mapsto \beta y$.

For $(\alpha_1,\dots,\alpha_n) \in \kxn$, there is an
automorphism of $\Obfqkn$ such that $x_i\mapsto \alpha_ix_i$
for all $i$.

For $(\alpha_1,\alpha_2,\beta_1,\beta_2) \in \kxfour$, there are
automorphisms of $\OqMtwo$ and $\OqGLtwo$ such that $X_{ij} \mapsto
\alpha_i\beta_j X_{ij}$ for all $i,j$. In other words,
$$\pmatrix X_{11} &X_{12} \\ X_{21}
&X_{22} \endpmatrix \longmapsto
\pmatrix \alpha_1 &0 \\ 0& \alpha_2
\endpmatrix \pmatrix X_{11} &X_{12} \\ X_{21}
&X_{22} \endpmatrix \pmatrix \beta_1 &0 \\ 0& \beta_2
\endpmatrix.$$

The automorphisms above do not all carry over to
$\OqSLtwo$ -- we must restrict attention to those which fix the
quantum determinant. For
$(\alpha,\beta)
\in \kxtwo$, there is an automorphism of
$\OqSLtwo$ such that 
$$\Xbar_{ij} \longmapsto \alpha^{3-2i}
\beta^{3-2j}\, \Xbar_{ij}$$
 for all $i,j$, that is,
$$\pmatrix \Xbar_{11} &\Xbar_{12} \\ \Xbar_{21}
&\Xbar_{22} \endpmatrix \longmapsto
\pmatrix \alpha &0 \\ 0& \alpha^{-1}
\endpmatrix \pmatrix \Xbar_{11} &\Xbar_{12} \\ \Xbar_{21}
&\Xbar_{22} \endpmatrix \pmatrix \beta &0 \\ 0& \beta^{-1}
\endpmatrix.$$
\enddefinition

\definition{1.4\. Example} The homogeneity of the plane in the
classical case carries over to its coordinate ring in the following
way -- if $M_1$ and $M_2$ are any maximal ideals of $k[x,y]$ of
codimension $1$ (these are the maximal ideals corresponding to points
in the plane with coordinates in $k$), there is an automorphism $\phi$
of $k[x,y]$ such that $\phi(M_1)= M_2$. Thus, if $k$ is algebraically
closed, the maximal ideals of $k[x,y]$ form a single orbit with
respect to the automorphisms of this algebra. 

While $\Oqktwo$ does not have enough automorphisms to map any maximal
ideal onto any other, there are still relatively large orbits.
Assuming that $k$ is algebraically closed and $q$ is not a root of
unity, the maximal ideals of $\Oqktwo$ can be seen in Example 1.2.
Using just the automorphisms defined in (1.3), there are three
orbits of maximal ideals:
$$\xalignat3 \bigl\{ \langle x,\, y-\beta\rangle &\bigm| \beta\in\kx
\bigr\} &\bigl\{ \langle x,\, y\rangle \bigr\}& 
& \bigl\{ \langle x-\alpha,\, y\rangle &\bigm| \alpha\in\kx
\bigr\}. \endxalignat$$
Note that each of the orbits above intersects to a prime ideal which
is stable under these automorphisms. The maximal ideals together with
these orbit-intersections account for all but one prime ideal of
$\Oqktwo$; for completeness, note that the remaining prime, namely
$\langle 0\rangle$, is also stable under the automorphisms.
\enddefinition

A similar pattern can be observed in $\OqSLtwo$, as follows.

\definition{1.5\. Example}
When $k$ is algebraically closed and $q$ is not a root of unity, the
prime spectrum of $\OqSLtwo$ can be displayed as shown below:

\ignore
$$\xymatrixrowsep{4pc}\xymatrixcolsep{2.25pc}
\xymatrix{
 &\ar@{--}[r] &\langle \Xbar_{11}-\lambda,\,\Xbar_{12},\, \Xbar_{21},\,
\Xbar_{22}-\lambda^{-1} \rangle \save+<0ex,-4ex> \drop{(\lambda \in
k^\times)} \restore \ar@{--}[r] & \\
 &&\langle \Xbar_{12},\,\Xbar_{21} \rangle \edge[ul] \edge[ur]\\
\langle \Xbar_{12} \rangle \edge[urr] &\edge[ur] \ar@{--}[r] &\langle
\Xbar_{12}-
\lambda \Xbar_{21} \rangle \save+<0ex,-4ex> \drop{(\lambda \in k^\times)}
\restore
\ar@{--}[r] &\edge[ul] &\langle \Xbar_{21}\rangle \edge[ull]\\
 &&\langle 0\rangle \edge[ull]\edge[ul] \edge[ur]\edge[urr]
}$$
\endignore

\noindent In this case, the maximal ideals form a single orbit under
the automorphisms described in (1.3), while the prime ideals of
height $1$ form three orbits. The two remaining primes can be
described as intersections of infinite orbits.
\enddefinition

Patterns analogous to those discussed in Examples 1.4 and 1.5 have
been found in all the other quantized coordinate rings introduced so
far, assuming that $k$ is algebraically closed and the parameters are
generic. With some modifications, the picture can be expanded to
include arbitrary infinite base fields. There is some interesting ring
theory which explains and predicts this behavior, and our main goal in
Part I is to present this theory. Some of the concepts used to
describe the picture only involve an arbitrary group of automorphisms
of a ring, but the key results hold when the group is an {\it algebraic
torus\/}, that is, a product of copies of the multiplicative group
$\kx$. The two examples above exhibit orbits of prime ideals which
intersect to stable prime ideals, which hints at the importance of
such orbit intersections. This hint leads to the key idea -- to group
prime ideals according to the intersections of their orbits with
respect to a specific group of automorphisms.

We begin with arbitrary actions of
groups on rings. Whenever we refer to a group acting on a ring, we shall
assume that it is acting by means of ring automorphisms (rather than
just by permutations or by invertible linear transformations,
for instance); similarly, actions on algebras are assumed to be
actions by algebra automorphisms. 

\definition{1.6\. $H$-prime ideals} Let $A$ be a ring, and let $H$ be a
group acting on $A$ (by automorphisms). Thus, we are given a
homomorphism $\phi:H\rightarrow \Aut A$, and we abbreviate
$\phi(h)(a)$ to $h(a)$ for $h\in H$ and $a\in A$. (Many authors write
$h.a$ for $\phi(h)(a)$.) For any ideal $P\vartriangleleft A$, set 
$$(P:H) \eqdef
\bigcap\limits_{h\in H} h(P),$$
the largest $H$-stable ideal of $A$ contained in $P$.

By restricting the usual definition of a prime ideal to $H$-stable
ideals, we obtain the concept of an {\it
$H$-prime\/} ideal of
$A$, namely any proper
$H$-stable ideal
$J$ of $A$ such that $I_1I_2 \not\subseteq J$ for all $H$-stable ideals
$I_1,\,I_2 \not\subseteq J$. In parallel with the notation $\spec A$,
we write $\Hspec A$ to denote the set of all $H$-prime ideals of $A$.
For example, if $q$ is not a root of unity, $k=\kbar$, and $H=\kxtwo$
acts as in (1.3), then
$$\Hspec \Oqktwo= \{\, \langle x,y\rangle,\, \langle
x\rangle,\,
\langle y\rangle,\, \langle 0\rangle\, \}.$$
\enddefinition

\proclaim{1.7\. Lemma} Let $H$ be a group acting on a ring $A$.

{\rm (a)} If $P$ is any prime ideal of $A$, then $(P:H)$ is an
$H$-prime ideal of $A$.

{\rm (b)} Now assume that $A$ is noetherian. Then a proper ideal
$J$ of $A$ is $H$-prime if and only if $J$ equals the intersection of
some finite $H$-orbit of prime ideals.

In particular, it follows that all $H$-prime ideals of $A$ are
semiprime in this case. \endproclaim

\demo{Proof} (a) Easy.

(b) E.g., see \cite{\BrGo, Lemma II.1.10}. \qed\enddemo

\definition{1.8\.  $H$-stratifications} Let $H$ be a group acting on
a ring
$A$. For each $H$-prime ideal $J$ of $A$, let
$$\specj A \eqdef \{ P\in \spec A \mid (P:H) =J\}.$$
This set is called {\it the $H$-stratum of
$\spec A$ corresponding to $J$\/}. In view of Lemma 1.7(a),
$$\spec A=
\bigsqcup\limits_{J\in\Hspec A} \specj A,$$
a partition that we call the {\it $H$-stratification\/} of $\spec A$.
\enddefinition

The $H$-stratifications just defined have similar properties to the
stratifications used in algebraic geometry, as follows.

\proclaim{1.9\. Lemma} Let $H$ be a group acting on a ring $A$.

{\rm (a)} The closure of each $H$-stratum in $\spec A$ is a union
of $H$-strata.

{\rm (b)} If $\Hspec A$ is finite, then each $H$-stratum is locally
closed in $\spec A$. \endproclaim

\demo{Proof} \cite{\GMurcia, Lemma 3.4}. \qed\enddemo

The stratification setup so far is extremely general, and we cannot
expect to prove much about it without specializing to cases with
additional hypotheses. One key specialization is to assume that $H$ is
an {\it affine algebraic group\/} over $k$, by which we just mean that
$H$ is isomorphic to a Zariski-closed subgroup of $GL_n(k)$ for some
$n$. Thus, $H$ is an affine algebraic variety as well as a group, and
the group operations are morphisms of varieties. We will not need much
at all of the general theory of algebraic groups, since we will
concentrate on one of the simplest kind, namely algebraic tori. To see
that a torus $(\kx)^r$ is
an algebraic group, note that it is isomorphic to the subgroup
of $GL_{r+1}(k)$ consisting of matrices $(a_{ij})$ satisfying the
equations $a_{ij}=0$ for $i\ne j$ and $a_{11}a_{22}
\cdots a_{r+1,r+1}= 1$.

\definition{1.10\. Rational actions} Let $A$ be a $k$-algebra, and $H$
a group acting on $A$. (As noted above, in this situation we assume
that $H$ acts on $A$ via $k$-algebra automorphisms.) Moreover, let us
assume that $H$ is an algebraic group over $k$. 

The action of $H$ on $A$ is said to be {\it rational\/} provided $A$
is a directed union of finite dimensional $H$-stable $k$-subspaces
$V_i$ such that the restriction maps $H\longrightarrow GL(V_i)$ are
morphisms (of algebraic groups), i.e., group homomorphisms which are
also morphisms of varieties. Fortunately for our purposes, the theory
of algebraic groups provides a nice criterion that allows us to see
quite easily when an action of a torus is rational, as follows.
\enddefinition

\definition{1.11\. Rational characters} Suppose that $H$ is an 
algebraic torus. Recall that a {\it character\/} of $H$ (with respect
to the base field
$k$) is any group homomorphism $H\rightarrow \kx$. Characters appear
whenever $H$ acts on a $k$-algebra $A$: If $x\in A$ is an
$H$-eigenvector (i.e., a simultaneous eigenvector for the actions of
all the automorphisms from $H$), then there is a character $\phi$ of
$H$ such that $h(x)= \phi(h)x$ for all $h\in H$. Of course, $\phi$ is
then called the {\it $H$-eigenvalue\/} of $x$.

 A character of $H$ is called
{\it rational\/} if it is also a morphism of varieties. Let
$X(H)$ denote the set of all rational characters of $H$; this is an
abelian group under pointwise multiplication, and it is easily
described. Namely, if
$H= (\kx)^r$, then $X(H)$ is a free abelian group in which the $r$
coordinate projections $(\kx)^r
\rightarrow \kx$ form a basis.
\enddefinition

\proclaim{1.12\. Theorem} Let $H$ be a torus acting on a $k$-algebra
$A$, and assume that $k$ is infinite. The action of $H$ on $A$ is
rational if and only if

{\rm (a)} The action is semisimple (i.e., $A$ is spanned by
$H$-eigenvectors); and

{\rm (b)} The $H$-eigenvalues for the $H$-eigenvectors in $A$ are
\underbar{rational} characters.
\endproclaim

\demo{Proof} \cite{\Nor, Chapter 5, Corollary to Theorem 36}.
\qed\enddemo

From a ring-theoretic point of view, conditions (a) and (b) of
Theorem 1.12 are the natural and useful conditions. Thus, we could
take them as our definition of a rational action of a torus, if
desired.

The next lemma illustrates one useful aspect of having a rational
action. (Recall that in general, an $H$-prime ideal in a noetherian
ring need only be semiprime.)

\proclaim{1.13\. Lemma} Suppose that $H$ is a torus, acting rationally
on a noetherian $k$-algebra $A$. Then every $H$-prime ideal of $A$
is prime. \endproclaim

\demo{Proof} If $J$ is an $H$-prime ideal of $A$, then Lemma 1.7(b)
implies that $A= (P:H)$ for some prime ideal $P$ whose $H$-orbit is
finite. Hence, the stabilizer subgroup $\stabh(P)$ has finite index in
$H$. Since $H$ acts rationally, the
map $H\rightarrow \spec A$ given by $h\mapsto h(P)$ is continuous with
respect to the Zariski topologies on $H$ and $\spec A$ \cite{\BrGo,
Lemma II.2.8}. Consequently, the set $V= \{ h\in H \mid h(P)\supseteq
P\}$ is closed in $H$. (We cannot say immediately that $\stabh(P)$ is
closed in $H$ because $\{P\}$ need not be closed in $\spec A$.) However,
because
$A$ is noetherian, any automorphism $\phi$ of $A$ for which $\phi(P)
\supseteq P$ must map $P$ onto itself. Hence, $V= \stabh(P)$, and thus
$\stabh(P)$ is indeed closed in $H$. 

Now $H$ is the disjoint union of the cosets of $\stabh(P)$. There are
only finitely many cosets, and they are all closed. However, as a
variety $H$ is irreducible (because its coordinate ring is a Laurent
polynomial ring, hence a domain), so it cannot be a finite union of
proper closed subsets. Thus $\stabh(P)=P$, that is, the $H$-orbit of
$P$ consists of $P$ alone. Therefore $J=P$, proving that $J$ is prime.
\qed\enddemo

We can now present a general theorem which provides a picture of the
structure of the $H$-stratification in our current setting. Recall
that a {\it regular\/} element in a noetherian ring is any
non-zero-divisor. We write $\Fract R$ to denote the Goldie quotient
ring of a semiprime noetherian ring $R$, and $Z(R)$ for the center of
a ring $R$.

\proclaim{1.14\. Stratification Theorem} \cite{\specstrat, \GMurcia}
Let
$A$ be a noetherian $k$-algebra, with
$k$ infinite, and let $H= (\kx)^r$ be a torus acting rationally on $A$.
For $J\in \Hspec A$, let $\E_J$ denote the set of all regular
$H$-eigenvectors in $A/J$.

{\rm (a)} $\E_J$ is a denominator set, and the localization $A_J \eqdef
(A/J)[\E_J^{-1}]$ is an $H$-simple ring (with respect to the induced
$H$-action).

{\rm (b)} $\spec_J A$ is homeomorphic to $\spec A_J$ via localization
and contraction, and $\spec A_J$ is homeomorphic to $\spec Z(A_J)$ via
contraction and extension.

{\rm (c)} $Z(A_J)$ is a Laurent polynomial ring of the form
$K_J[z_1^{\pm1},\dots,z_{n(J)}^{\pm1}]$, with $n(J)\le r$, over the
fixed field $K_J \eqdef Z(A_J)^H= Z(\Fract A/J)^H$.
\endproclaim

\demo{Proof} \cite{\BrGo, Chapter II.3}. \qed\enddemo

Of course, the theorem above does not say much if the $H$-strata are
very small and there are many of them. For instance, in the extreme
case the $H$-strata might be singletons, in which case the theorem is
trivial. To get the most information out of this picture, we would
like there to be only finitely many $H$-strata, so that the
$H$-stratification breaks up the prime spectrum into relatively large
sets. Many quantized coordinate rings are iterated skew polynomial
extensions of $k$, and the following theorem can be applied to those
algebras.

\proclaim{1.15\. Theorem} \cite{\specstrat, \BrGo} Let $A$ be an
iterated skew polynomial algebra 
$$k[x_1] [x_2;\tau_2,\delta_2] \cdots
[x_n;\tau_n,\delta_n],$$
 and let $H$ be a group acting on
$A$, such that $x_1,\dots,x_n$ are $H$-eigenvectors.
Assume that there exist $h_1,\dots,h_n \in H$ such that:

{\rm (a)} $h_i(x_j)= \tau_i(x_j)$ for $i>j$; and

{\rm (b)} The $h_i$-eigenvalue of $x_i$ is not a root of unity
for any $i$.

\noindent Then $A$ has at most $2^n$ $H$-prime ideals. Moreover, if
$H$ is a torus acting rationally on $A$, then for each $J\in\Hspec
A$, the field $K_J$ (from part {\rm (c)} of the Stratification
Theorem) equals
$k$.
\endproclaim

\demo{Proof} \cite{\BrGo, Theorems II.5.12 and II.6.4}. \qed\enddemo

One of the tools involved in proving the Stratification Theorem is an
equivalence between rational $(\kx)^r$-actions and $\ZZ^r$-gradings,
part of which we now sketch. Since this is intended to be applied to
the $H$-prime factor algebras $A/J$, we work with an algebra called
$B$ rather than $A$.

\definition{1.16\. Actions versus gradings} Suppose that $B$ is a
noetherian
$k$-algebra, with $k$ infinite, and that a torus $H=(\kx)^r$
acts rationally on $B$. Because of Theorem 1.12,
$$B= \bigoplus_{g\in X(H)} B_g,$$
where $B_g$ denotes the $H$-eigenspace of $B$ with eigenvalue $g$.
Since $H$ acts by automorphisms, $B_gB_{g'} \subseteq B_{gg'}$ for all
$g,g'\in G$, that is, $B$ is graded by the group $X(H)\cong \ZZ^r$.
(Conversely, any grading of a $k$-algebra by $\ZZ^r$ corresponds to a
rational action of $(\kx)^r$ on the algebra.) Problems concerning the
$H$-action translate into problems concerning the grading in the
following way:
$$\align H\text{-eigenvectors}\quad
&\longleftrightarrow\quad \text{homogeneous elements}\\
H\text{-stable ideals}\quad &\longleftrightarrow\quad
\text{homogeneous ideals}\\
H\text{-prime ideals}\quad &\longleftrightarrow\quad
\text{graded-prime ideals}. \endalign$$

To prove part (a) of the Stratification Theorem, we need to be able to
localize an $H$-prime ring $B$ with respect to its regular
$H$-eigenvectors and obtain an $H$-simple ring. Translating to the
graded case, we need to localize a graded-prime ring with respect to
its homogeneous regular elements and obtain a graded-simple ring. In
other words, what is required is a version of Goldie's Theorem for the
setting of graded rings. This cannot be obtained in general -- there
are easy examples of commutative, noetherian, semiprime $\ZZ$-graded
rings where the localization with respect to all homogeneous regular
elements is not graded-simple. For our present purposes, it suffices
to consider prime graded rings, for which the following theorem is
available.
\enddefinition

\proclaim{1.17\. Graded Goldie Theorem} Let $G$ be an
abelian group, and let $R$ be a $G$-graded, graded-prime, right
graded-Goldie ring. Let $\E$ be the set of all homogeneous regular
elements in $R$.
Then $\E$ is a right denominator set, and $R[\E^{-1}]$ is a
graded-simple, graded-artinian ring. \endproclaim

\demo{Proof} \cite{\GoSt, Theorem 1}. \qed\enddemo

Theorem 1.17 moves us to the setting of graded-simple rings, and the
prime ideals in such rings can be analyzed as follows.

\proclaim{1.18\. Proposition} Let $G$ be an abelian
group, and let $R$ be a $G$-graded, graded-simple ring.

{\rm (a)} $\spec R$ is homeomorphic to $\spec Z(R)$ via contraction and
extension.

{\rm (b)} If $G\cong \ZZ^r$, then $Z(R)$ is a Laurent polynomial ring,
in at most $r$ indeterminates, over the field $Z(R)_1$ (the identity
component of $Z(R)$).
\endproclaim

\demo{Proof} \cite{\BrGo, Lemma II.3.7 and Proposition II.3.8}.
\qed\enddemo

Let us conclude Part I by presenting an open problem.

\definition{1.19\. Problem} Suppose that $A$ is a noetherian
$k$-algebra, and that a torus
$H=(\kx)^r$ acts rationally on $A$. Find conditions which imply that
$A$ has only finitely many $H$-primes. These conditions should be
\roster
\item"$\bullet$"  Reasonably easy to verify; and
\item"$\bullet$" Satisfied by all the standard examples.
\endroster
In other words, we would like to have a theorem which we can apply to
quantized coordinate rings without masses of long calculations. In
seeking such a theorem, a warning is in order: When the parameters are
roots of unity, quantized coordinate rings usually have
infinitely many
$H$-primes. Thus, whatever hypotheses might be used in a solution to
this problem will have to correspond to the generic situation when
applied to quantized coordinate rings.
\enddefinition

\head II. Primitive Ideals \endhead

We now concentrate on primitive ideals as opposed to general prime
ideals, and on ways to organize the {\it primitive spectrum\/} of a
ring $A$. This set, denoted $\prim A$, is the set of all
left primitive ideals of $A$. We view $\prim A$ as a
topological space equipped with the Zariski topology, so that $\prim
A$ is a subspace of $\spec A$.

The question whether the left primitive ideals and the right primitive
ideals coincide in a noetherian ring remains open. To avoid this
problem, we shall use the term {\it primitive ideal\/} to refer only
to {\it left\/} primitive ideals.

In a classical coordinate ring over an algebraically closed field, the
maximal ideals correspond to points of the underlying variety. A naive
geometric analogy in the noncommutative world would be to view the
maximal ideals in a ring as points of a `noncommutative variety'.
However, experience in ring theory teaches us that there are
too few maximal ideals in general to hold sufficient information.
Further, the influence of representation theory leads us to study the
primitive ideals, as one key to irreducible representations (i.e.,
simple modules). Let us consider our simplest example, $\Oqktwo$. 

\definition{2.1\. Example} Assume that
$k$ is algebraically closed and $q$ is not a root of unity. The prime
ideals of $\Oqktwo$ are displayed in Example 1.2. Observe that the
only maximal ideals are the ideals
$$\langle x-\alpha,\,y\rangle \qquad\quad \text{and} \qquad\quad
\langle x,\,y-\beta\rangle,$$
for $\alpha,\beta\in k$. Comparing these with the maximal ideals in
the classical coordinate ring $\Oktwo$, we see that the maximal ideals
of $\Oqktwo$ correspond only to points on the $x$- and $y$-axes of
$k^2$. From this point of view, the remainder of the $xy$-plane has
been `lost'. 

As suggested above, let us widen our view to include all the primitive
ideals. In $\Oqktwo$, there is one non-maximal primitive ideal, namely
$\langle 0\rangle$. Thus, comparing $k^2$ with $\prim\Oqktwo$, we can
now say that the points on the $x$- and $y$-axes correspond precisely
to the maximal ideals of $\Oqktwo$, while all other points of $k^2$
correspond (not bijectively, of course) to the zero ideal. We could
say that the off-axis part of $k^2$ has `collapsed' to a single point.
Later, we shall  elaborate this point of view further.
\enddefinition

Since the primitive ideals of an algebra $A$ are (by definition) the
annihilators of the simple $A$-modules, it would seem that to
determine these primitive ideals, we should find all the simple
$A$-modules and then calculate their annihilators. However, to find all
the simple modules over an infinite dimensional algebra is usually an
impossible problem. As a substitute, Dixmier promulgated the following
program for enveloping algebras of Lie algebras: Find the primitive
ideals first, and then for each primitive ideal, find at least one
simple module having that annihilator. In order to carry out this
program, we must be able to detect the primitive ideals without knowing
the simple modules in advance, and so some criterion other then the
definition is required. Since all primitive ideals are prime, the
question becomes, how can we tell which prime ideals are primitive?
Dixmier developed two criteria, one of which is purely algebraic and one
of which is phrased topologically, as follows.

\definition{2.2\. Rational and locally closed primes} Let $P$ be a
prime ideal in a noetherian
$k$-algebra $A$. First, we say that $P$ is {\it rational\/} if and
only if $Z(\Fract A/P)$ is algebraic
over $k$. Secondly, we say that $P$ is {\it locally closed\/} provided
$P$ is a locally closed point in $\spec A$, i.e., the singleton set
$\{P\}$ is closed in some neighborhood of $P$. This condition may be
rephrased as follows: $P$ is locally closed if and only if
$$\bigcap \{Q\in \spec A \mid Q
\supsetneq P\} \supsetneq P.$$
Thus, $P$ is locally closed if and only if in the prime ring $R/P$, the
intersection of all nonzero prime ideals is nonzero. Prime rings with the
latter property are sometimes called {\it G-rings\/}, in which case
locally closed primes are called {\it G-ideals\/}.
\enddefinition

\proclaim{2.3\. Theorem} Let
$\gfrak$ be a finite dimensional Lie algebra over a field of
characteristic zero. Then
$$\align \prim U(\gfrak) &= \{ \text{\rm locally closed prime ideals
of\ } U(\gfrak)\} \\
 &= \{\text{\rm rational prime ideals of\ } U(\gfrak)\} \endalign$$
\endproclaim

\demo{Proof} This theorem was originally proved by Dixmier
\cite{\Dix} and Moeglin \cite{\Moe} assuming an
algebraically closed base field. Their result was extended to
the non-algebraically closed case by Irving and Small \cite{\IrSm}.
\qed\enddemo

\definition{2.4\. The Dixmier-Moeglin equivalence} We say that an
algebra $A$ satisfies the {\it Dixmier-Moeglin equivalence\/}
if the conclusion of Theorem 2.3 holds in $A$, that is, the primitive,
locally closed, and rational prime ideals of $A$ all coincide.
\enddefinition

There are some relations among these three types of prime ideals which
hold under fairly general hypotheses. One such hypothesis is the
following adaptation of Hilbert's Nullstellensatz to noncommutative
noetherian algebras.

\definition{2.5\. The noncommutative Nullstellensatz} A 
$k$-algebra
$A$ is said to satisfy the {\it Null\-stel\-len\-satz over
$k$\/} if and only if
\roster
\item"(a)" The Jacobson radical of every factor ring of $A$ is nil; and
\item"(b)" $\End_A(M)$ is
algebraic over $k$ for all simple $A$-modules $M$.
\endroster
If $A$ is noetherian, condition (a) is equivalent to $A$
being a Jacobson ring, i.e., $J(A/P)=0$ for all $P\in
\spec A$. \enddefinition

The Nullstellensatz is essentially automatic if the field is large
enough. In particular:

\proclaim{2.6\. Proposition} \cite{\Ami} If $k$ is uncountable,
then every countably generated $k$-algebra satisfies the
Nullstellensatz over $k$. \endproclaim

\demo{Proof} \cite{\McRo, Corollary 9.1.8}.
\qed\enddemo

For algebras over countable fields, the following theorem is often
useful.

\proclaim{2.7\. Theorem} Suppose that a $k$-algebra
$A$ has subalgebras
$A_0=k \subseteq A_1\subseteq \cdots\subseteq A_t=A$ such that for all
$i>0$, either $A_i$ is a finitely generated
$A_{i-1}$-module on each side, or $A_i$
is a homomorphic image of a skew polynomial ring
$A_{i-1}[x_i;\tau_i,\delta_i]$.  Then $A$ satisfies the
Nullstellensatz over $k$. \endproclaim

\demo{Proof} This is a special case of \cite{\McRo, Theorem 9.4.21}.
\qed\enddemo

For prime ideals in a noetherian algebra satisfying the
Nullstellensatz, the following general implications are known:
$$\text{locally closed} \implies \text{primitive} \implies
\text{rational}$$
\cite{\BrGo, Lemma II.7.15}.
Closing the loop (i.e., proving that `rational$\implies$locally
closed') is usually the most difficult part of establishing that an
algebra satisfies the Dixmier-Moeglin equivalence. In the situation of
the Stratification Theorem, it is advantageous to bring the torus action
into the loop -- this helps in the proofs, and supplies an additional
criterion for primitivity, namely the condition that a prime ideal be
a maximal element of its $H$-stratum. Then, two implications need to
be proved to close the loop, namely
$$\text{rational} \implies \text{maximal in stratum} \implies
\text{locally closed},$$
but the second is quite easy. The precise theorem is as follows.

\proclaim{2.8\. Theorem} \cite{\specstrat} Let $A$ be a
noetherian
$k$-al\-gebra with $k$ infinite, and let $H=(\kx)^r$ be a torus
acting rationally on $A$. Assume that $\Hspec A$ is finite, and that 
$A$ satisfies the Nullstellensatz over $k$. Then
$$\align \prim A &= \{\text{\rm locally closed prime ideals of\ } A\}
\\
 &= \{\text{\rm rational prime ideals of\ } A\} \\
 &= \bigsqcup_{J\in\Hspec A} \{\text{\rm maximal elements of\ } \specj
A\}. \endalign$$
Moreover, if $k$ is algebraically closed, the $H$-orbits in $\prim A$
coincide with the
$H$-strata $\prim_J A \eqdef (\prim A)\cap (\specj A)$.
\endproclaim

\demo{Proof} \cite{\BrGo, Theorems II.8.4 and II.8.14}. \qed\enddemo

Of the three criteria for primitivity given in this theorem, the third
is typically easiest to apply. Here is an illustration.

\definition{2.9\. Example} Let us return to our second basic
example,
$\OqSLtwo$, assuming that $k$ is algebraically closed and $q$ is not
a root of unity. Recall from (1.3) the rational action of $H=\kxtwo$
on this algebra. The prime spectrum of
$\OqSLtwo$ was displayed in Example 1.5, and the action of $H$ on
these prime ideals is easy to determine. In particular, there are
only four $H$-prime ideals in $\OqSLtwo$, namely $\langle
\Xbar_{12},\,\Xbar_{21} \rangle$, $\langle \Xbar_{12}\rangle$,
$\langle \Xbar_{21}\rangle$, and $\langle 0\rangle$. Thus, there are
four $H$-strata in
$\spec\OqSLtwo$, which we display as follows.

\ignore
$$\xymatrixrowsep{4pc}\xymatrixcolsep{2.25pc}
\xymatrix{
 &\ar@{--}[r] &\langle \Xbar_{11}-\lambda,\,\Xbar_{12},\, \Xbar_{21},\,
\Xbar_{22}-\lambda^{-1} \rangle \save+<0ex,-4ex> \drop{(\lambda \in
k^\times)} \restore \ar@{--}[r] & \\
 &&\langle \Xbar_{12},\,\Xbar_{21} \rangle \edge[ul] \edge[ur]\\
\langle \Xbar_{12}\rangle & \ar@{--}[r] &\langle
\Xbar_{12}-
\lambda \Xbar_{21} \rangle \save+<0ex,-4ex> \drop{(\lambda \in k^\times)}
\restore
\ar@{--}[r] &&\langle \Xbar_{21}\rangle\\
 &&\langle 0\rangle \edge[ul] \edge[ur]
}$$
\endignore

\bigskip

The Nullstellensatz holds for
$\OqSLtwo$ by Theorem 2.7. Hence, Theorem 2.8 im\-plies that
$\prim\OqSLtwo$ consists of all prime ideals except for $\langle
\Xbar_{12},\Xbar_{21}\rangle$ and $\langle 0\rangle$. Moreover,
there are precisely four $H$-orbits in $\prim\OqSLtwo$:
$$\matrix &\bigl\{ \langle \Xbar_{11}-\lambda,\,\Xbar_{12},\,
\Xbar_{21},\,
\Xbar_{22}-\lambda^{-1} \rangle \bigm| \lambda\in\kx \bigr\} \\  \\
\bigl\{ \langle \Xbar_{12} \rangle \bigr\} &\bigl\{ \langle
\Xbar_{12}- \lambda \Xbar_{21} \rangle \bigm| \lambda\in\kx \bigr\}
&\bigl\{ \langle \Xbar_{21} \rangle \bigr\} \endmatrix$$
\enddefinition

Now that we have access to finding the primitive ideals in quantized
coordinate rings, let us turn to the global problem -- trying to
understand the primitive spectrum of such an algebra $A$ as a whole. We
would like $\prim A$ to reflect some kind of `noncommutative
geometry'. Since there is as yet no indication of what might play the
role of regular functions on $\prim A$, we focus for now on the
topological structure of this space.
\medskip

{\sl For the remainder of Part II, assume that $k$ is
algebraically closed.}

\definition{2.10\. Problem} Let $V$ be an affine variety over $k$, with
classical coordinate ring $\O(V)$, and suppose that
$A$ is some quantized coordinate ring of $V$. Since $\max \O(V)
\approx V$, we may view $\prim A$ as a `quantization of $V$'. Then
the problem arises: How are $\prim A$ and $V$ related?
\enddefinition

\definition{2.11\. Example} Assume that $q$ is not a root of unity. In
Example 2.1, we suggested that
$\prim\Oqktwo$ could be viewed as the union of the $x$-axis, the
$y$-axis, and one other point obtained from collapsing the rest of the
$xy$-plane. This leads to a map $\phi$ from
$k^2$ onto $\prim\Oqktwo$, as in the following sketch:

\ignore
$$\xymatrixrowsep{0.1pc}\xymatrixcolsep{0.1pc}
\xymatrix{
 &&&&& &&&&& &&&&& &&&&&\place \edge[8,-4] \\ 
 \\  \\  \\
 &&&&\place \edge[8,-4] \edge[-4,16] \\ 
  \\  \\  \\
 &&&&& &&&&& &&&&& &\place \\
  \\  \\  \\
\place \edge[-4,16] \\
  \\
 &&&&& &&&&&\place\ar@{->>}[6,0]_{\phi} \\
 \\  \\  \\  \\  \\
 &&&&& &&&&&\place \\
 \\  \\
 &&\place \edge[4,16] &&&&& &&&&& &&&&& &\place \\  \\
 &&&&& &&&&&\bullet \\  \\
 &&\place \edge[-4,16] &&&&& &&&&& &&&&& &\place \\  \\
 &&&&& &&&&&\bullet
}$$
\endignore

\bigskip

\noindent To describe this map more precisely, recall that
$\prim\Oqktwo$ consists of the maximal ideals
$\langle x-\alpha,\,y\rangle$ and $\langle x,\,y-\beta\rangle$, for
$\alpha,\beta\in k$, together with $\langle0\rangle$.
Thus, $\phi$ is given as follows:
$$\align (\alpha,0) &\longmapsto \langle x-\alpha,\,y\rangle \\
(0,\beta) &\longmapsto \langle x,\,y-\beta\rangle \\
\text{other\ } (\alpha,\beta) &\longmapsto \langle0\rangle\ .
\endalign$$ 
It is easy to check that $\phi$ is continuous. In fact, the topology
on $\prim\Oqktwo$ equals the quotient topology induced by
$\phi$. Thus, $\prim\Oqktwo$ is a topological quotient of $k^2$.
\enddefinition

\proclaim{2.12\. Conjecture} If an algebra $A$ is one of the
`standard' quantized coordinate rings of an affine variety $V$, then
$\prim A$ is a topological quotient of $V$. \endproclaim

This conjecture is known to hold in several cases:

\quad(1)  $A=\OqSLtwo$, when $q$ is not a root of unity. We invite the
reader to try this as an exercise.

\quad(2) $A=\Obfqkxn \eqdef \Obfqkn[x_1^{-1},\dots,x_n^{-1}]$,
for arbitrary $\bfq$. This follows from work of De Concini-Kac-Procesi
\cite{\DKP}, Hodges \cite{\Hod}, Vancliff \cite{\Van}, Brown-Goodearl
\cite{\BroGoo}, Goodearl-Letzter \cite{\qaffine}, and others.

\quad(3) $A=\Obfqkn$ and a more general type of algebra known as a
`quantum toric variety' (which we will describe below), when the
subgroup of $\kx$ generated by the entries of $\bfq$ does not contain
$-1$. This is work of Goodearl and Letzter \cite{\qaffquo}.

\quad(4) $A=\O_q(\frak{sp}\,k^4)$, the single-parameter quantized
coordinate ring of symplectic 4-space, when $q$ is not a root of
unity. This algebra was first defined in \cite{\RTF}; for somewhat
simpler presentations see \cite{\Mus} and \cite{\Hortpaper}. The
topological quotient here was established by Horton \cite{\Hortdiss,
Theorem 7.9}.

\definition{2.13\. Quantum tori} Given a multiplicatively antisymmetric
$n\times n$ matrix $\bfq= (q_{ij})$ over $k$, the corresponding {\it
quantized coordinate ring of $\kxn$\/} is the $k$-algebra
$$\Obfqkxn \eqdef k\langle x_1^{\pm1},\dots,x_n^{\pm1} \mid x_ix_j=
q_{ij}x_jx_i
\text{\ for all\ } i,j\rangle.$$
The torus $H=(\kx)^n$ acts rationally on the algebra $A=\Obfqkxn$ in
the same way as it acts on $\Obfqkn$. As is easily checked,
$\langle0\rangle$ is the only $H$-prime in $A$, and $\prim A$ is a
single $H$-stratum as well as a single $H$-orbit. Hence, for any
primitive ideal $P$, there is a bijection
$$H/\stabh(P) \longleftrightarrow \prim A.$$
This bijection is a homeomorphism, assuming that $H/\stabh(P)$ is
given the quotient topology. Thus, the fact that $\prim A$ is a
topological quotient of $H$ is easily established in this case.
\enddefinition

\definition{2.14\. Quantum affine spaces} Now let $A=\Obfqkn$,
and recall that
$H=(\kx)^n$ acts rationally on $A$ by $k$-algebra automorphisms such
that
$(\alpha_1,\dots,\alpha_n).x_i=
\alpha_ix_i$ for $(\alpha_1,\dots,\alpha_n) \in H$ and $i=1,\dots,n$.
Let $W$ be the collection of subsets of $\{1,\dots,n\}$. There is a
bijection
$$\align W &\longrightarrow \Hspec A \\
w &\longmapsto J_w \eqdef \langle x_i\mid i\in w\rangle. \endalign$$
Thus, the $H$-stratifications of $\spec A$ and $\prim A$, and the
localizations of $A$ appearing in the Stratification Theorem, are
indexed by the $H$-primes
$J_w$. To simplify notation, we re-index using $W$. In particular, we
write
$$\align \prim_wA &\eqdef \prim_{J_w} A= \{P\in \prim A\mid x_i\in P
\iff i\in w\}\\
A_w &\eqdef A_{J_w} = (A/J_w)[x_j^{-1}\mid j\notin w] \endalign$$
for $w\in W$. Note that each $A_w$ is a quantum torus.

The torus $H$ acts on $\Okn$ exactly as it does on $A$; this action is
induced from the action of $H$ on $k^n$ by the rule
$$(\alpha_1,\dots,\alpha_n).(a_1,\dots,a_n) \eqdef (\alpha_1^{-1}a_1,
\dots, \alpha_n^{-1}a_n).$$
There are $2^n$ $H$-orbits in $k^n$, which we index by $W$ as follows:
$$(k^n)_w \eqdef \{ (a_1,\dots,a_n)\in k^n \mid a_i=0
\iff i\in w \}$$
for $w\in W$. We may note that $(k^n)_w$ is isomorphic to a torus of
rank $n-|w|$. The results discussed in (2.13) imply that
$\prim_w A$ is a topological quotient of $(k^n)_w$ for each $w$. Thus,
the problem is to patch individual topological quotient maps $(k^n)_w
\twoheadrightarrow \prim_w A$ together, to obtain a topological
quotient map $k^n
\twoheadrightarrow \prim A$. Our solution to this problem requires a
small technical condition, phrased in terms of $\langle
q_{ij}\rangle$, the subgroup of $\kx$ generated by the entries
$q_{ij}$ of $\bfq$.
\enddefinition

\proclaim{2.15\. Theorem} Assume that either
$-1\notin
\langle q_{ij}\rangle$ or $\chr k =2$. Then there exist compatible, 
$H$-equivariant topological quotient maps
$$k^n \twoheadrightarrow \prim \Obfqkn \qquad\qquad \text{and}
\qquad\qquad \spec \O(k^n) \twoheadrightarrow \spec \Obfqkn$$
such that for  $w\in W$, the inverse image of $\prim_w\Obfqkn$
is $(k^n)_w$. Moreover, the fibres over points in $\prim_w\Obfqkn$ are
$G_w$-orbits in
$(k^n)_w$ for certain subgroups $G_w\subseteq H$.
\endproclaim

\demo{Proof} \cite{\qaffquo, Theorem 4.11; \GDurham, Theorem 3.5}.
These papers also describe how to calculate the subgroups $G_w$.
\qed\enddemo

To illustrate this theorem, we use a single parameter quantum affine
$3$-space, this being the simplest case in which the topological
quotient map differs from what one might naively write down.

\definition{2.16\. Example} Choose a non-root of unity
$q\in\kx$, and let $A=\O_q(k^3)$. Then the entries $q_{ij}$ in $\bfq$
consist of
$q$,
$q^{-1}$, and $1$, and so the group $\langle q_{ij}\rangle$ is
infinite cyclic. In particular, $-1\notin
\langle q_{ij}\rangle$ unless $\chr k =2$. Now let $p$ be one of
the square roots of $q$ in $\kx$. The topological quotient map
$k^3\twoheadrightarrow
\prim A$ given by Theorem 2.15 can be described as shown below, where
all $\lambda_i\in\kx$.
$$\gather (0,0,0) \longmapsto \langle x_1,\, x_2,\, x_3 \rangle \\
\alignedat2  (\lambda_1,0,0) &\longmapsto \langle x_1-\lambda_1,\,
x_2,\, x_3 \rangle &\qquad (\lambda_1,\lambda_2,0) &\longmapsto
\langle x_3 \rangle\\ 
(0,\lambda_2,0) &\longmapsto \langle x_1,\, x_2-\lambda_2,\, x_3
\rangle &\qquad (\lambda_1,0,\lambda_3) &\longmapsto \langle x_2
\rangle\\ 
(0,0,\lambda_3) &\longmapsto \langle x_1,\, x_2,\, x_3-\lambda_3
\rangle &\qquad (0,\lambda_2,\lambda_3) &\longmapsto \langle x_1
\rangle  \endalignedat \\
(\lambda_1,\lambda_2,\lambda_3) \longmapsto \langle \lambda_2x_1x_3-
p\lambda_1\lambda_3x_2 \rangle\ . \endgather$$
Note the appearance of $p$ in the final line -- without that factor,
the resulting map from $k^3$ to $\prim A$ will still be surjective,
but not Zariski-continuous.
\enddefinition

\definition{2.17} We indicate one basic mechanism from the proof of
Theorem 2.15. Given a multiplicatively antisymmetric $n\times n$
matrix
$\bfq$, we write parallel, coordinate-free descriptions of the
algebras 
$R= \O(k^n)$ and $A= \Obfqkn$ as follows. Namely, $R=
k(\ZZ^+)^n$, a semigroup algebra, and $A=k^c(\ZZ^+)^n$, a twisted
semigroup algebra for a suitable cocycle $c: \ZZ^n\times
\ZZ^n\longrightarrow \kx$. There are many choices of $c$; we just need
to have $c(\epsilon_i,\epsilon_j)c(\epsilon_j,\epsilon_i)^{-1}
=q_{ij}$ for all $i,j$, where $\epsilon_1,\dots,\epsilon_n$ is the
standard basis for $\ZZ^n$. Now both $R$ and $A$ have bases identified
with $(\ZZ^+)^n$, and so there is a \underbar{vector space}
isomorphism $\Phi_c: A\longrightarrow R$ which is the identity on
$(\ZZ^+)^n$. Similarly, $\Phi_c$ extends to a vector space isomorphism
from the group algebra $k\ZZ^n$ onto the twisted group algebra
$k^c\ZZ^n$, and so for each $w\in W$ we obtain a vector space
isomorphism $\Phi_c$ from $A_w$ onto a subalgebra $R_w$ of $k\ZZ^n$. The
key to Theorem 2.15 is to choose $c$ so that the above maps behave well:
\enddefinition

\proclaim{2.18\. Key Lemma} There is a choice of cocycle $c$ such that
$\Phi_c$ yields
\underbar{$k$-algebra} maps $Z(A_w) \longrightarrow R_w$ for all $w$.

For this choice of $c$, the topological quotient maps $\max R
\twoheadrightarrow \prim A$ and $\spec
R\twoheadrightarrow
\spec A$ can be described by the rule
$$Q \longmapsto \bigl( \text{\rm the largest ideal of $A$ contained in\
} \Phi_c^{-1}(Q) \bigr).$$
\endproclaim

\demo{Proof} The first statement follows from \cite{\qaffquo, (4.2),
(4.6--4.8), (3.5)}, while the second is \cite{\GDurham, Lemma 3.6}.
\qed\enddemo

For a more precise description of this map in terms of
operations within the commutative algebra $R$, see \cite{\qaffquo,
\GDurham}.

Since the method just sketched is based on twisting the
polynomial ring $k(\ZZ^+)^n$ by a cocycle, it readily extends to a
somewhat more general class of  algebras twisted by cocycles.

\definition{2.19\. Cocycle twists} Suppose that $G$ is a group, and
that $R$ is a
$G$-graded $k$-algebra. Let $c:G\times
G\rightarrow
\kx$ be a $2$-cocycle, normalized so that $c(1,1)=1$ (or $c(0,0)=1$,
in case $G$ is written additively). The {\it twist of $R$ by $c$\/} is
a $k$-algebra based on the same $G$-graded vector space as $R$, but
with a new multiplication $*$ defined on homogeneous elements as
follows: $r*s \eqdef c(\alpha,\beta)rs$ for $r\in R_\alpha$ and $s\in
R_\beta$. 

Now specialize to the case that $R$ is a commutative affine $G$-graded
$k$-algebra, and $A$ is the twist of $R$ by a 2-cocycle $c$. Then $R$
is generated by finitely many homogeneous elements, say
$y_1,\dots,y_n$, of degrees
$\alpha_1,\dots,\alpha_n$. The algebra $A$ is generated by the same
elements
$y_1,\dots,y_n$, and $y_i*y_j= q_{ij}y_j*y_i$ for all $i,j$, where
$q_{ij}= c(\alpha_i,\alpha_j)c(\alpha_j,\alpha_i)^{-1}$. Consequently,
$A\cong \Obfqkn/I$ for $\bfq=(q_{ij})$ and some ideal $I$.

In particular, if $G=\ZZ^n$ and $\dim R_\alpha =1$ for all
$\alpha\in G$, then
$R$ is the coordinate ring of an affine toric variety $V$, and we
regard $A$ as a quantized coordinate ring of $V$. This case was
studied by Ingalls \cite{\Ing}, who introduced the term {\it quantum
toric variety\/} to describe the resulting algebras $A$.

The constructions behind Theorem 2.15 adapt well to factor algebras
$\Obfqkn/I$, and that theorem extends to the cocyle twisted setting as
follows.
\enddefinition

\proclaim{2.20\. Theorem} Let $G$ be a
torsionfree abelian group, and let
$R$ be a commutative, affine, $G$-graded
$k$-algebra. Let $A$ be the
twist of $R$ by a 2-cocycle $c:G\times
G\rightarrow \kx$. Assume that $-1\notin \langle
\text{\rm image}(c)\rangle \subseteq \kx$, or that $\chr k=2$. 

Then there
exist compatible topological quotient maps
$$\max R \twoheadrightarrow \prim A \qquad\qquad \text{and}
\qquad\qquad \spec R \twoheadrightarrow \spec A,$$
which are equivariant with respect to the action of a suitable torus.
\endproclaim

\demo{Proof} \cite{\qaffquo, Theorem 6.3; \GDurham, Theorem 4.5}.
\qed\enddemo

It is not clear whether the hypothesis concerning $-1$ can be removed
from Theorems 2.15 and 2.20. We end Part II by putting an extreme case
forward as an open problem.

\definition{2.21\. Problem} Assume that $\chr k\ne 2$. Consider the
single parameter algebras 
$$\O_{-1}(k^n)= k\langle x_1,\dots,x_n \mid
x_ix_j= -x_jx_i \text{\ for all\ } i\ne j \rangle.$$
The methods used to prove Theorem 2.15 still work for $\O_{-1}(k^2)$
and $\O_{-1}(k^3)$. These methods break down for $\O_{-1}(k^4)$, but
extensive ad hoc calculations lead to a Zariski-continuous surjection
$k^4 \twoheadrightarrow \prim \O_{-1}(k^4)$; in higher dimensions,
the problem is completely open. Thus, we ask:

For $n\ge 4$, is the space $\prim \O_{-1}(k^n)$ a topological quotient
of $k^n$?
\enddefinition

\head III. Quantum Matrices \endhead

The focus on topological quotients in Part II was chosen to emphasize
one way in which the quantized coordinate ring of a variety can be
geometrically similar to the classical coordinate ring. We can also
ask about algebraic similarities, of which there are many -- chain
conditions, homological conditions, etc. In fact, there exist much
tighter similarities -- many classical theorems have surprisingly
close quantum analogs, once they are properly rephrased. We
illustrate this principle by discussing quantum matrices, that is,
the quantized coordinate rings of varieties of matrices. The
$2\times2$ case was presented in (1.1); we now give the general
definition.

\definition{3.1\. Generators and relations} Let $n$ be a positive
integer and $q\in\kx$. The {\it quantized coordinate ring of $n\times n$
matrices\/} with parameter $q$ is the $k$-algebra with generators
$X_{ij}$ for
$i,j=1,\dots,n$ such that for all $i<l$ and $j<m$, the generators
$X_{ij},X_{im},X_{lj},X_{lm}$ satisfy the defining relations for
$\OqMtwo$. As in (1.1), five of these relations can be summarized in
the following mnemonic diagram:
\ignore $$\xymatrix{
X_{ij} \ar[r]^q \ar[d]_q &X_{im} \ar[d]^q \ar@{<->}[dl]_1  \\
X_{lj}\ar[r]^q &X_{lm}
}$$ \endignore
The remaining relation is $X_{ij}X_{lm}-X_{lm}X_{ij}=
(q-q^{-1})X_{im}X_{lj}$.

The $n\times n$ quantum determinant is modelled on the usual
determinant, but with powers of $-1$ replaced by powers of $-q$. More
precisely, the {\it $n\times n$ quantum determinant\/} is the
element 
$$D_q \eqdef \sum_{\pi\in S_n} (-q)^{\ell(\pi)} X_{1,\pi(1)}
X_{2,\pi(2)}
\cdots X_{n,\pi(n)}\ \ \in\ \ \OqMn,$$
where $S_n$ denotes the symmetric group and $\ell(\pi)$, the {\it
length\/} of a permutation $\pi$, is the minimum length for an
expression of $\pi$ as a product of simple transpositions
$(i,i{+}1)$. It is known that $D_q$ lies in the center of $\OqMn$.
Hence, one defines quantized coordinate rings $\OqGLn \eqdef
\OqMn[D_q^{-1}]$ and $\OqSLn \eqdef \OqMn/ \langle D_q-1 \rangle$ as
before.

The algebra $\OqMn$ is a bialgebra with comultiplication and counit
maps
$$\Delta: \OqMn \longrightarrow \OqMn\otimes\OqMn
\qquad\text{and}\qquad \varepsilon:
\OqMn \longrightarrow k$$
 such that  $\Delta(X_{ij}) = \sum_{l=1}^n X_{il}\otimes X_{lj}$
and $\varepsilon(X_{ij}) = \delta_{ij}$
for all $i,j$.
In particular, $\Delta(D_q)= D_q\otimes D_q$ and $\varepsilon(D_q)=1$.
\enddefinition

All this structure is exactly parallel to the classical case, which
we get if $q=1$. Much interesting geometry has resulted from viewing
sets of matrices of a given size as algebraic varieties and focusing
on constructs from linear algebra as geometric processes. One such
line leads to determinantal ideals, as follows.

\definition{3.2\. Classical determinantal ideals} Let $t\le n$ be
positive integers, and consider the variety
$$V_t \eqdef \{ n\times n \text{\ matrices of rank\ } <t \},$$
the closed subvariety of the affine space $M_n(k)$ defined by the
vanishing of all $t\times t$ minors. From linear algebra, $V_t$ is
the image of the matrix multiplication map
$$M_{n,t-1}(k)\times M_{t-1,n}(k) \longrightarrow M_n(k).$$
Since  $M_{n,t-1}(k)$ and $M_{t-1,n}(k)$ are irreducible varieties, it
follows that $V_t$ is irreducible.

Let $I_t \vartriangleleft \O(M_n(k))$ be the ideal of polynomial
functions vanishing on $V_t$, so that $\O(M_n(k))/I_t= \O(V_t)$.
On the geometric side, $V_t$ is defined by the vanishing of all
$t\times t$ minors. However, this only tells us that $I_t$ equals the
\underbar{radical} of the ideal generated by the $t\times t$ minors.
It is a classical theorem that these minors actually generate this
ideal:
\enddefinition

\proclaim{3.3\. Theorem} $I_t$ equals the
ideal of
$\O(M_n(t))$ generated by all $t\times t$ minors. \endproclaim

\demo{Proof} See, e.g., \cite{\BrVe, \DeCEP}. \qed\enddemo

\proclaim{3.4\. Corollary} The set of all $t\times t$ minors in
$\O(M_n(k))$ generates a prime ideal. \qed\endproclaim

In the quantum world, there is no variety $V_t$, and so we cannot ask
for a direct analog of Theorem 3.3. However, there are analogs of
minors, which means that we can look for an analog of Corollary 3.4.

\definition{3.5\. Quantum minors} Let $I,J\subseteq \{1,\dots,n\}$ be
index sets with
$|I|=|J|=t$. We may write the elements of these sets in ascending
order, say $I= \{i_1<\cdots<i_t\}$ and $J= \{j_1<\cdots<j_t\}$ for
short. There is a natural $k$-algebra embedding
$\phi_{I,J} :
\O_q(M_t(k))
\rightarrow \OqMn$ such that $\phi_{I,J}(X_{lm})= X_{i_lj_m}$ for all
$l,m$. The {\it quantum minor with index sets $I$ and $J$\/} is the
element
$$[I|J] \eqdef \phi_{I,J}(D_q^{t\times t})\ \ \in\ \ \OqMn,$$
where $D_q^{t\times t}$ denotes the quantum determinant in
$\O_q(M_t(k))$.
\enddefinition

\proclaim{3.6\. Theorem} The ideal $I_t$
of $\OqMn$ generated by all $t\times t$ quantum minors is completely
prime, i.e.,
$\OqMn/I_t$ is an integral domain. \endproclaim

\demo{Proof} \cite{\GLenDuke, Theorem 2.5}. \qed\enddemo

Although many steps in the proof of the classical result have no
analogs in the quantum case, one part of the classical pattern does
carry over, as we now summarize.

\definition{3.7} As noted above,
$V_t$ is the image of the multiplication map
$$\mu:
M_{n,t-1}(k)\times M_{t-1,n}(k) \longrightarrow M_n(k).$$
Hence, the ideal $I_t$ is the kernel of the comorphism 
$$\mu^*:
\O(M_n) \longrightarrow
\O(M_{n,t-1}\times M_{t-1,n}).$$
We may identify $\O(M_{n,t-1}\times M_{t-1,n})$ with $\O(M_{n,t-1})
\otimes \O(M_{t-1,n})$, which allows us to describe $\mu^*$ as the
composition of the maps
$$\O(M_n) @>{\,\,\Delta\,\,}>> \O(M_n)\otimes\O(M_n)
@>{\,\,\text{quo}\otimes\text{quo}\,\,}>> \O(M_{n,t-1}) \otimes
\O(M_{t-1,n}).$$
\enddefinition

\definition{3.8\. A quantum analog} Quantized coordinate rings for the
rectangular matrix varieties $M_{n,t-1}(k)$ and $M_{t-1,n}(k)$ may be
defined as the subalgebras of $\OqMn$ generated by those $X_{ij}$ with
$j<t$ (respectively, $i<t$). There are natural $k$-algebra retractions
of $\OqMn$ onto these subalgebras, and so
$$\align \O_q(M_{n,t-1}(k)) &\cong \OqMn/ \langle X_{ij} \mid j\ge
t\rangle \\
\O_q(M_{t-1,n}(k)) &\cong \OqMn/ \langle X_{ij} \mid i\ge
t\rangle. \endalign$$
Thus, the quantum analog of the comorphism $\mu^*$ in (3.7) is the
$k$-algebra map
$$\mu_q^* \eqdef \Oq(M_n) @>{\,\,\Delta\,\,}>> \Oq(M_n)\otimes\Oq(M_n)
@>{\,\,\text{quo}\otimes\text{quo}\,\,}>> \O_q(M_{n,t-1}) \otimes
\O_q(M_{t-1,n}).$$
It is easy to check that $\O_q(M_{n,t-1}(k)) \otimes
\O_q(M_{t-1,n}(k))$ is an iterated skew polynomial algebra over $k$,
and therefore a domain. Thus, to prove that the ideal $I_t$ of
$\OqMn$ is completely prime, one just has to show that
$I_t=\ker(\mu_q^*)$. This is the heart of the proof of Theorem 3.6.
\enddefinition

To understand quantum analogs of other geometric aspects of matrices,
and also to understand the quantum matrix algebra better, we would
like to know its prime and primitive ideals. We approach this problem
via the Stratification Theorem, as discussed in Parts I and II.

{\sl For the remainder of Part III, assume that
$q$ is not a root of unity, and set
$A=\OqMn$.}

\definition{3.9\. Problem}  In parallel with the $2\times2$ case discussed in (1.3),
the torus
$H=(\kx)^{2n}$ acts on
$A$ by $k$-algebra automorphisms so that
$$(\alpha_1,\dots,\alpha_n,\beta_1,\dots,\beta_n).X_{ij}=
\alpha_i\beta_j X_{ij}$$
 for all $i,j$. These automorphisms are called `winding
automorphisms', because they arise from the bialgebra structure on
$A$ in a manner analogous to the classical winding automorphisms on
enveloping algebras of Lie algebras.
According to Theorem 1.15, there are at most $2^{n^2}$
$H$-prime ideals in $A$. The basic problem is:

Determine the $H$-primes of $A$.
\enddefinition
 
\definition{3.10\. Example} The $2\times2$ case of Problem 3.9 is
easily solved -- there are exactly 14 $H$-prime ideals in
$\OqMtwo$, as displayed in the following diagram. Each
$2\times2$ pattern here is shorthand for a set of generators of an
ideal -- a bullet $\bullet$ in position $(i,j)$ corresponds to a
generator $X_{ij}$; a circle $\circ$ in a given position is a
placeholder; and the square $\square$ denotes the $2\times2$ quantum
determinant.

\bigskip

\centerline{\epsfbox{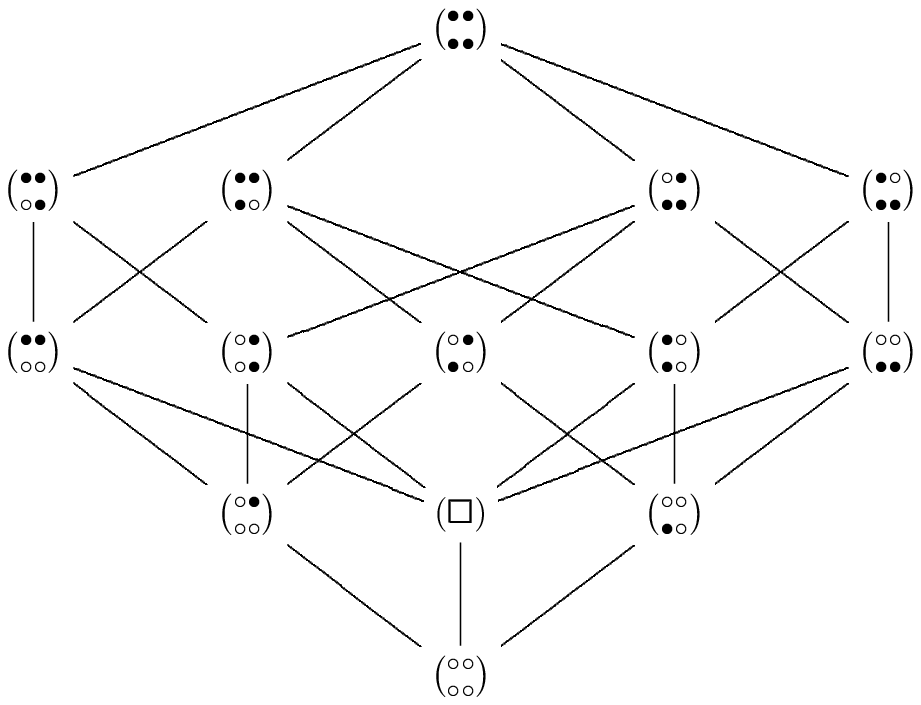}}
\enddefinition

\definition{3.11} Certain types of $H$-primes in $A$ are already
known. For convenient labelling, we carry over the term `rank' from
ordinary matrices to the quantum case, as follows.  We define the
{\it rank\/} of a prime ideal $P$ in $A$ to be the minimum $r$ such
that $P$ contains all $(r+1)\times(r+1)$ quantum minors.

The $H$-primes of $A$ of rank $n$ are those which do not contain the
quantum determinant
$D_q$. By localization, these correspond to the $H$-primes of $\OqGLn$,
and it is known that those, in turn, correspond to the $H$-primes of
$\OqSLn$. The latter can be determined using results of Hodges and
Levasseur
\cite{\HoLethree, \HoLeN}. In particular, $\OqSLn$ has $(n!)^2$
$H$-primes, parametrized by
$S_n\times S_n$, and it follows from work of Joseph \cite{\Jsur,
Th\'eor\`eme 3} that each of these
$H$-primes is generated by a set of quantum minors. We conclude that
back in $A$, the
$H$-primes of rank $n$ are generated -- up to localization at the
powers of $D_q$ -- by sets of quantum minors. 

At the other extreme, the $H$-primes of rank at most 1 are the
$H$-primes of
$A$ which contain all
$2\times 2$ quantum minors. These were determined by Goodearl and
Lenagan \cite{\GLenMurcia, Proposition 3.4}. There are $(2^n-1)^2+1$
such $H$-primes, all having the form
$$\bigl \langle [I|J] \bigm| |I|=|J|=2 \bigr \rangle
+\bigl \langle X_{ij}\mid i\in R \bigr \rangle +\bigl \langle
X_{ij}\mid j\in C \bigr \rangle$$ 
for $R,C\subseteq \{1,\dots,n\}$. For these
$H$-primes, we have generating sets consisting of quantum minors,
since each
$X_{ij}$ is a $1\times1$ quantum minor.
\enddefinition

\proclaim{3.12\. Conjecture} Every $H$-prime of $A$ is generated
by a set of quantum minors. \endproclaim

It is easily seen that the conjecture holds when $n=2$, in view of
(3.10). Cauchon proved that there are enough quantum
minors to separate the $H$-primes in $A$: For any $H$-primes $P
\subsetneq Q$, there is a quantum minor in $Q\setminus P$ 
\cite{\Cau, Proposition 6.2.2 and Th\'eor\`eme 6.3.1}. The
$3\times3$ case of the conjecture has been established by Goodearl
and Lenagan \cite{\GLenJalg, Theorem 7.4}, and the
$n\times n$ case, assuming that $k=\CC$ and $q$ is transcendental over
$\QQ$, has been proved by Launois \cite{\Lau, Th\'eor\`eme 3.7.2}. We
shall display the solution to the $3\times3$ case below.

Cauchon's and Launois's results are existence theorems -- they do not
provide descriptions of which sets of quantum minors generate
$H$-primes. Such descriptions are needed not only for completeness,
but also to get full benefit from the stratification, e.g., to
determine the prime and primitive ideals in each $H$-stratum via
Theorems 1.14 and 2.8. Thus, we accompany Conjecture 3.12 with the
following problems.

\definition{3.13\. Problems} If $J\in \Hspec A$, then Theorems 1.14
and 1.15 tell us that the center of the localization
$A_J= (A/J)[\E_J^{-1}]$ is a Laurent polynomial ring
$k[z_1^{\pm1},\dots,z_{n(J)}^{\pm1}]$.

(a) Assuming $J$ is generated by a set $\M$ of quantum minors, find a
formula for $n(J)$ in terms of $\M$.

(b) Find explicit descriptions of the indeterminates $z_i\in A_J$.
\enddefinition

We now summarize the solution to the $3\times3$ case of Conjecture
3.12 given in \cite{\GLenJalg}.

\definition{3.14} As in (3.11), we may divide up the $H$-primes in
$\O_q(M_3(k))$ according to their ranks. Those of ranks 0, 1, and 3
were known earlier, while the ones of rank 2 were first determined in
\cite{\GLenJalg}. The numerical count is as follows:
$$\alignat2 \rank\ 0 &:\quad &1& \\
\rank\ 1 &: \quad &49& \\
\rank\ 2 &: \quad &144& \\
\rank\ 3 &: \quad &36& \\
\text{total\hphantom{ra}} &: \quad &230&.  \endalignat$$
The determination of these $H$-primes was done partly by {\it ad
hoc\/} methods, which are unlikely to work in the general case. In
particular, Cauchon has given a formula for the total number of
$H$-primes in
$\OqMn$ \cite{\Cau, Th\'eor\`eme 3.2.2 and Proposition 3.3.2}, which
shows that
$\O_q(M_4(k))$ has 6902
$H$-primes!
\enddefinition

\definition{3.15} The $H$-primes in $\O_q(M_3(k))$ can be displayed as
in the following diagrams, where each $3\times3$ pattern represents a
set of generators for an $H$-prime. As in (3.10), circles are
placeholders and bullets represent generators $X_{ij}$. This time,
squares and rectangles represent $2\times2$ quantum minors whose row
and column index sets correspond to the edges. Finally, the diamond
that appears in four patterns represents the
$3\times3$ quantum determinant. Below are samples showing the ideals
corresponding to two patterns.
\ignore
$$\align \vcenter{ \xymatrixrowsep{0.15pc}\xymatrixcolsep{0.15pc}
\xymatrix{ \cbb \\  \hhrzvrt &*{\phcirc}+0 \\  \hhrz
&*{\phcirc}+0 }}
\quad &\longleftrightarrow\quad\ 
\langle\, X_{12},\, X_{13},\, [23|12],\, [23|13],\, [23|23]\, \rangle
\\
 \vcenter{ \xymatrixrowsep{0.15pc}\xymatrixcolsep{0.15pc}
\xymatrix{ \plc &\place \edge[dl] \edge[dr] &\plc \\  \place
\edge[dr] &\plc &\place \edge[dl] &*{\phcirc}+0 \\  \plb &\place
&\plc  }}
\quad &\longleftrightarrow\quad\ 
\langle\, [123|123],\, X_{31}\, \rangle. \endalign$$
\endignore
\enddefinition

The case of rank $0$ is trivial -- there is only one $H$-prime of
this rank, corresponding to the following pattern:

\ignore
$$\xymatrixrowsep{0.15pc}\xymatrixcolsep{0.15pc}
\xymatrix{ \bbb \\  \bbb \\  \bbb }$$
\endignore

\medskip

The 49 $H$-primes of rank $1$ correspond to the patterns in Figure A
below. 
\goodbreak
\midinsert
\centerline{\epsfbox{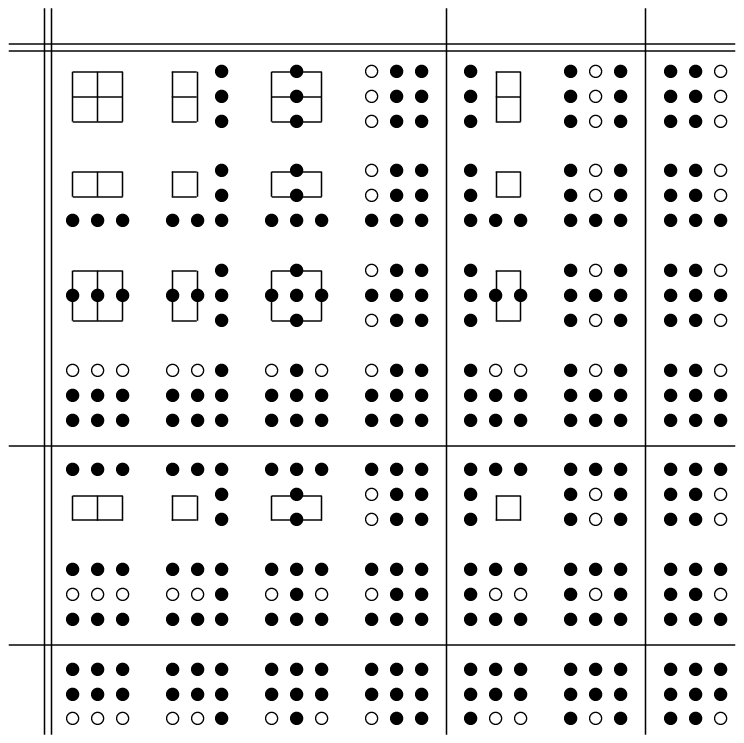}}
\smallskip
\centerline{\bf Figure A}
\endinsert

As indicated in (3.11), the $H$-primes of maximal rank were
known, up to localization, from the results of
\cite{\HoLethree, \HoLeN, \Jsur}. The sets of quantum minors
which generate the corresponding $H$-primes in $\O_q(SL_3(k))$ also,
as it turns out, generate the $H$-primes of rank $3$ in
$\O_q(M_3(k))$. These 36 ideals correspond to the patterns in Figure
B.
\goodbreak
\midinsert
\centerline{\epsfbox{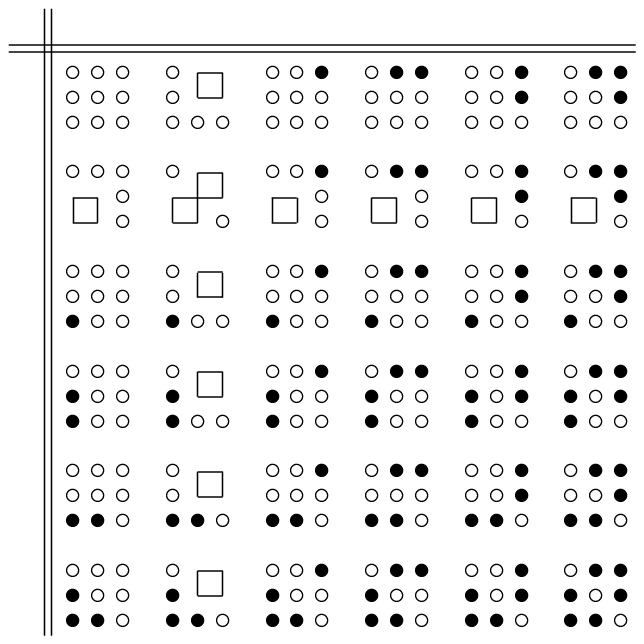}}
\smallskip
\centerline{\bf Figure B}
\endinsert

The 144 patterns for generating sets of the
$H$-primes of rank $2$ in $\O_q(M_3(k))$ are given in our
final display, Figure C.
The procedure used in \cite{\GLenJalg} to determine these $H$-primes
involved three steps. First, some general theory developed in
\cite{\GLenIJM} provided a reduction mechanism relating the
$H$-primes in $\OqMn$ to pairs of $H$-primes from smaller quantum
matrix algebras. Consequently, we could find the $H$-primes in
$\O_q(M_3(k))$ from the (known) $H$-primes in $\OqMtwo$, but only as
kernels of certain algebra homomorphisms. This process also gave
precise counts for the number of $H$-primes of each rank. In the
second step, the information from Step 1 was used to determine the
quantum minors contained in each
$H$-prime, thus yielding at least potential sets of generators.
Finally, the third step consisted of proving that each set of quantum
minors appearing in Step 2 does generate an $H$-prime, and that the
resulting $H$-primes are distinct. Since that yielded a list of the
correct number of
$H$-primes, we were done.

It may be useful to break the general problem down into steps of
similar type, as follows.

\definition{3.16\. Problems} We return to the general $n\times n$
situation, keeping the field $k$ arbitrary but still requiring $q$ to
be a non-root of unity.

 (a) Which sets $\M$ of quantum minors in
$\OqMn$ generate prime ideals?

(b) Develop general theorems to prove that suitable ideals of the form
$\langle\M\rangle$ are prime.

(c) Find combinatorial data to parametrize the sets $\M$ in (a).
\enddefinition

We conclude by stating a result which illustrates one pattern which
solutions to the above problems might take. Part (a) is an easy
exercise involving the relations among the $X_{ij}$, part (b)
can be proved by showing that $\OqMn/\langle\X\rangle$ is an iterated
skew polynomial extension of $k$, and part (c) is another easy
exercise.

\proclaim{3.17\. Sample result} {\rm (a)} If $P$ is an $H$-prime ideal
of $\OqMn$, then the set $\X= P\cap \{X_{ij}\mid i,j=1,\dots,n\}$
satisfies the following condition:
\roster
\item"(*)" If $X_{ij}\in \X$, then either $X_{lm}\in\X$ for all $l\ge
i$ and $m\le j$, or else $X_{lm}\in\X$ for all $l\le
i$ and $m\ge j$.
\endroster

{\rm (b)} If $\X$ is any subset of $\{X_{ij}\mid i,j=1,\dots,n\}$ which
satisfies {\rm (*)}, then $\X$ generates an $H$-prime ideal of $\OqMn$,
and $\langle\X\rangle \cap \{X_{ij}\mid i,j=1,\dots,n\} =\X$.

{\rm (c)} Given subsets $I,J \subseteq \{1,\dots,n\}$ and
nondecreasing functions $ f: \{1,\dots,n\}\setminus J \rightarrow
\{2,\dots,n+1\}\setminus I$ and $g: \{1,\dots,n\}\setminus I
\rightarrow \{2,\dots,n+1\}\setminus J$, the set 
$$\align \X(I,J,f,g) \eqdef\ &\bigl\{ X_{ij} \bigm| i\in I \bigr\}
\cup \bigl\{ X_{ij} \bigm| i\notin I;\ j\notin J;\ i\ge f(j)
\bigr\} \\ 
&\cup \bigl\{ X_{ij} \bigm| j\in J \bigr\} \cup \bigl\{ X_{ij}
\bigm| i\notin I;\ j\notin J;\ j\ge g(i) \bigr\} \endalign$$
satisfies {\rm (*)}.
Conversely, any subset $\X\subseteq \{X_{ij}\mid i,j=1,\dots,n\}$
which satisfies {\rm (*)} equals $\X(I,J,f,g)$ for some $I$, $J$,
$f$, $g$.
\qed\endproclaim

\goodbreak
\topinsert
\centerline{\epsfbox{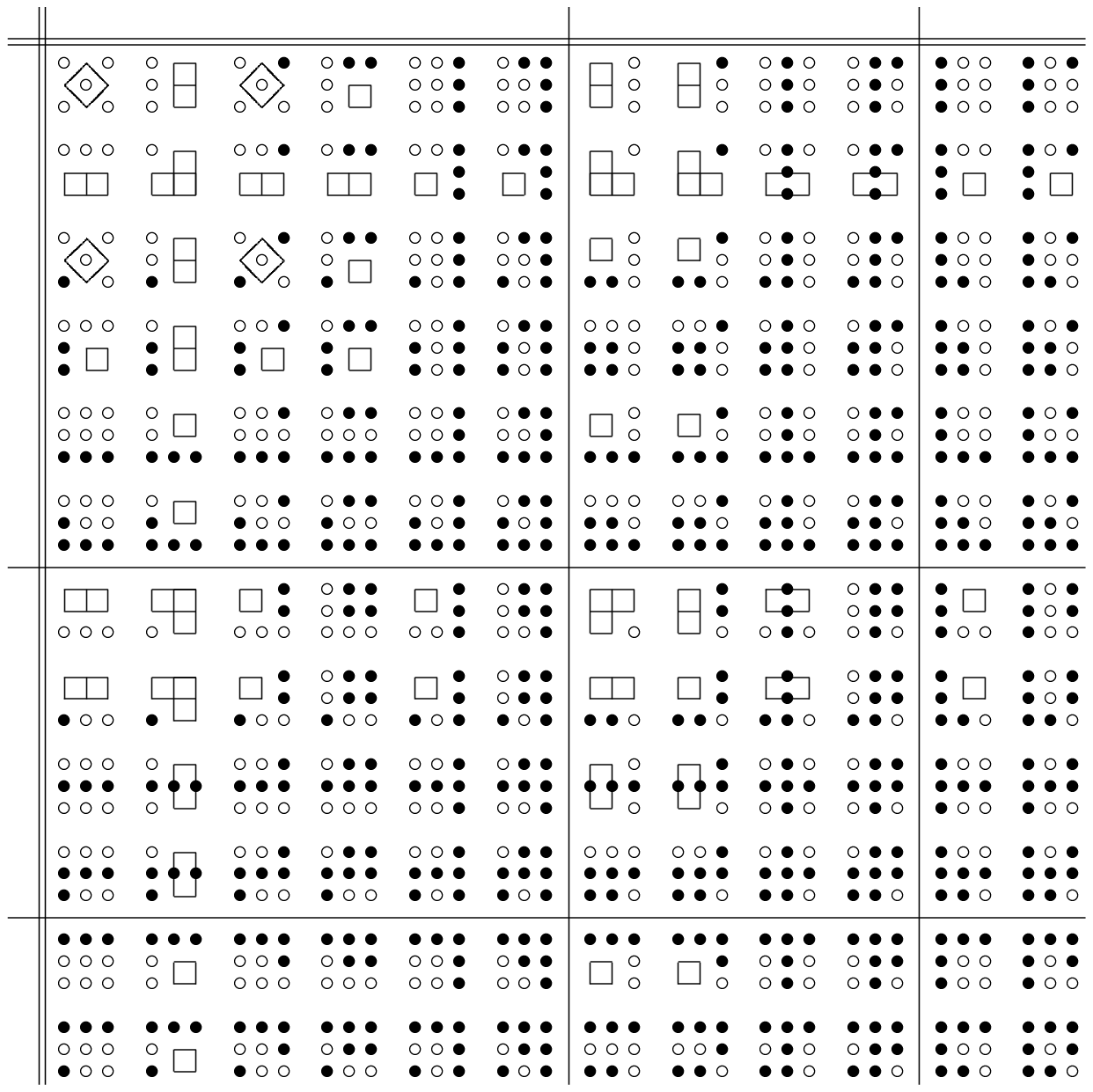}}
\smallskip
\centerline{\bf Figure C}
\endinsert


\Refs\widestnumber\no{{\bf 99}}

\ref\no\Ami \by  A. S. Amitsur \paper Algebras over infinite fields
\jour Proc. Amer. Math. Soc. \vol 7 \yr 1956 \pages 35-48 \endref

\ref\no \BroGoo \by K. A. Brown and K. R. Goodearl \paper Prime
spectra of quantum semisimple groups \jour
Trans\. Amer\. Math\. Soc\. \vol 348 \yr 1996 \pages 2465--2502
\endref

\ref\no\BrGo \bysame \book Lectures on
Algebraic Quantum Groups \bookinfo Advanced Courses in Mathematics CRM
Bar\-ce\-lo\-na \publaddr Basel \yr 2002 \publ Birkh\"auser \endref

\ref\no\BrVe \by W. Bruns and U. Vetter \book Determinantal Rings
\bookinfo Lecture Notes in Math. 1327
\publaddr Berlin
\yr 1988 \publ Springer-Verlag
\endref

\ref\no\Cau \by G. Cauchon \paper Spectre premier de $\OqMn$. Image
canonique et s\'eparation normale \jour J. Algebra \toappear \endref

\ref\no\ChPr \by V. Chari and A. Pressley \book A Guide to Quantum
Groups \publaddr Cambridge \yr 1994 \publ Cambridge Univ. Press
\endref

\ref\no\DeCEP \by C. De Concini, D. Eisenbud, and C. Procesi \paper
Young diagrams and determinantal varieties \jour Invent. Math. \vol 56
\yr 1980 \pages 129-165\endref 

\ref\no\DKP \by C. De Concini, V. Kac, and C. Procesi \paper
Some remarkable degenerations of quantum groups \jour
Comm\. Math\. Phys\. \vol 157 \yr 1993 \pages 405--427 \endref

\ref\no\Dix \by J. Dixmier \paper Id\'eaux primitifs dans les
alg\`ebres enveloppantes \jour J. Algebra \vol 48 \yr 1977 \pages
96--112 \endref

\ref\no\GMurcia \by K. R. Goodearl \paper Prime spectra of quantized
coordinate rings \inbook in Interactions between Ring Theory and
Representations of Algebras (Murcia 1998) \eds F.
Van Oystaeyen and M. Saor\'\i n \publaddr New York \yr 2000
\publ Dekker \pages 205-237
\endref

\ref\no\GDurham \by \bysame \paper Quantized primitive ideal spaces as
quotients of affine algebraic varieties \inbook in Quantum Groups and
Lie Theory \ed A. Pressley
\bookinfo London Math. Soc. Lecture Note Series 290 \pages 130-148
\publaddr Cambridge \yr 2001
\publ Cambridge Univ. Press
\endref

\ref\no\GLenMurcia \by K. R. Goodearl and T. H. Lenagan \paper Prime
ideals in certain quantum determinantal rings \inbook in Interactions
between Ring Theory and Representations of Algebras (Murcia 1998) \eds
F. Van Oystaeyen and M. Saor\'\i n \publaddr New York \yr 2000
\publ Dekker \pages 239-251
\endref

\ref\no\GLenDuke \bysame \paper Quantum determinantal ideals \jour
Duke Math. J. \vol 103 \yr 2000 \pages 165-190
\endref

\ref\no\GLenIJM \bysame \paper Prime ideals invariant under winding
automorphisms in quantum matrices \jour Internat. J. Math. \vol 13
\yr 2002 \pages 497-532
\endref

\ref\no\GLenJalg \bysame \paper Winding-invariant prime ideals in
quantum
$3\times3$ matrices \jour J. Algebra \toappear \endref

\ref\no\qaffine \by K. R. Goodearl and E. S. Letzter \paper Prime and
primitive spectra of multiparameter quantum affine spaces \paperinfo
in Trends in Ring Theory (Miskolc, 1996) (V. Dlab and L. Marki, eds.)
\jour Canad. Math. Soc. Conf. Proc. Series \vol 22
\yr 1998 \pages 39-58\endref

\ref\no\specstrat \by \bysame \paper The Dixmier-Moeglin equivalence
in quantum coordinate rings and quantized Weyl algebras \jour Trans.
Amer. Math. Soc. \vol 352 \yr 2000 \pages 1381-1403 \endref

\ref\no\qaffquo \by \bysame \paper Quantum $n$-space as a quotient of
classical $n$-space \jour Trans. Amer. Math. Soc. \vol 352 \yr 2000
\pages 5855-5876 \endref

\ref\no\GoSt \by K. R. Goodearl and J. T. Stafford \paper The
graded version of Goldie's Theorem \paperinfo in Algebra and its
Applications (Athens, Ohio, 1999) (D. V. Huynh, S. K. Jain, and S. R.
L\'opez-Permouth, Eds.) \jour Contemp. Math. \vol 259
\yr 2000 \pages 237-240
\endref

\ref\no\Hod \by T. J. Hodges \paper Quantum tori and Poisson tori
\paperinfo Unpublished notes (1994) \endref

\ref\no\HoLethree \by T. J. Hodges and T. Levasseur\paper Primitive
ideals of
${\bold C}_q[SL(3)]$ \jour Comm. Math. Phys. \vol 156 \yr 1993
\pages 581-605
\endref

\ref\no\HoLeN \by \bysame \paper Primitive
ideals of
${\bold C}_q[SL(n)]$ \jour J. Algebra \vol 168 \yr 1994 \pages 455-468
\endref

\ref\no\Hortdiss \by K. L. Horton \paper Prime spectra of iterated
skew polynomial rings of quantized coordinate type \finalinfo Ph.D.
Dissertation (2002) University of California at Santa Barbara \endref

\ref\no\Hortpaper \bysame \paper The prime and primitive spectra of
multiparameter quantum symplectic and Euclidean spaces \jour Communic. in
Algebra \toappear
\endref

\ref\no\Ing \by C. Ingalls \paper Quantum toric varieties \toappear
\endref

\ref\no\IrSm \by R. S. Irving and L. W. Small \paper On the
characterization of primitive ideals in enveloping algebras \jour
Math. Z. \vol 173 \yr 1980 \pages 217--221 \endref

\ref\no\Jsur \by A. Joseph \paper Sur les ideaux g\'en\'eriques de
l'alg\`ebre des fonctions sur un groupe quantique \jour C. R. Acad.
Sci. Paris, S\'er. I
\vol 321 \yr 1995 \pages 135-140
\endref

\ref\no\Kas \by C. Kassel \book Quantum Groups \bookinfo Grad.
Texts in Math. 155
\publaddr New York \yr 1995 \publ Springer-Verlag \endref

\ref\no\KlSc \by A. Klimyk and K. Schm\"udgen \book Quantum
Groups and their Representations \publaddr Berlin \yr 1997 \publ
Springer-Verlag
\endref

\ref\no\Lau \by S. Launois \paper Les id\'eaux premiers invariants de
$O_q({\Cal M}_{m,p}(\Bbb C))$ \finalinfo Preprint (2002) \endref

\ref\no\McRo \by J. C. McConnell and J. C. Robson \book
Noncommutative Noetherian Rings \publ Wiley-Interscience \publaddr
Chi\-ches\-ter-New York \yr 1987 \moreref \bookinfo Reprinted with
corrections, Grad. Studies in Math. 30 \publaddr Providence \publ
Amer. Math. Soc. \yr 2001
\endref

\ref\no\Moe \by C. Moeglin \paper Id\'eaux primitifs des alg\`ebres
enveloppantes \jour J. Math. Pures Appl. \vol 59 \yr 1980 \pages
265--336 \endref

\ref\no\Mus \by I. M. Musson \paper Ring theoretic properties of the
coordinate rings of quantum symplectic and Euclidean space\inbook in
Ring Theory, Proc. Biennial Ohio State--Denison Conf., 1992\eds S. K.
Jain and S. T. Rizvi\bookinfo \publaddr Singapore\yr 1993\publ World
Scientific\pages 248-258\endref

\ref\no\Nor \by D. G. Northcott \book Affine Sets and Affine Groups
\bookinfo London Math. Soc. Lecture Note Series 39 \publaddr
Cambridge
\yr 1980 \publ Cambridge Univ. Press \endref

\ref\no\RTF \by N. Yu. Reshetikhin, L. A. Takhtadjan, and L. D.
Faddeev\paper Quantization of Lie groups and Lie algebras\jour
Leningrad Math. J.\vol 1\yr 1990\pages 193-225 \endref

\ref\no\Van \by M. Vancliff \paper Primitive and Poisson spectra of
twists of polynomial rings \jour Algebras and Representation Theory
\vol 2 \yr 1999 \pages 269-285
\endref


\endRefs
\enddocument